\documentclass[11pt,reqno,a4paper]{amsart}
\usepackage{color}
\usepackage{amssymb,amsmath}
\usepackage[latin1]{inputenc}
\usepackage[active]{srcltx}
\usepackage{graphicx}
\usepackage{textgreek}
\usepackage{epsfig,graphics}
\usepackage{psfrag}
\usepackage{caption}
\usepackage{subcaption}
\def\N{\mathbb{N}}
\def\R{\mathbb{R}}
\def\m1{{I\!\!M}}

\newcommand{\ci}{\mathbb{C}}
\newcommand{\mb}[1]{\mathbb{#1}}

\renewcommand{\to}{\rightarrow}
\newcommand{\pa}{\partial}

\newcommand{\ino}{\int_{\Omega}}

\newcommand{\rife}[1]{(\ref{#1})}
\newcommand{\ov}[1]{\overline{#1}}
\newcommand{\un}[1]{\underline{#1}}

\newcommand{\scp}{\scriptstyle}
\newcommand{\sscp}{\scriptscriptstyle}
\newcommand{\dsp}{\displaystyle}

\renewcommand{\dfrac}{\displaystyle\frac}
\newcommand{\finedim}{\hspace{\fill}$\square$}
\newcommand{\intbar}{\mathop{\int\makebox(-15.5,0){\rule[6pt]{.7em}{0.3pt}}%
\kern-6pt}\nolimits}

\newcommand{\ii}{\infty}

\newcommand{\eps}{\varepsilon}
\newcommand{\dt}{\delta}

\newcommand{\al}{\alpha}

\newcommand{\sg}{\sigma}

\newcommand{\om}{\Omega}
\newcommand{\lm}{\lambda}

\newcommand{\ssb}{\sscp \om}

\newcommand{\prl}{{\textbf{(}\mathbf P_{\mathbf \lm}\textbf{)}}}
\newcommand{\fq}{f_{\scp q}}

\renewcommand{\rho}{\mbox{\Large \textrho}}

\newtheorem{theorem}{Theorem}[section]
\newtheorem{proposition}[theorem]{Proposition}
\newtheorem{lemma}[theorem]{Lemma}
\newtheorem{corollary}[theorem]{Corollary}
\newtheorem{remark}[theorem]{Remark}
\newtheorem{definition}[theorem]{Definition}
\newcommand{\brm}{\begin{remark}\rm}
\newcommand{\erm}{\end{remark}}
\newcommand{\bdf}{\begin{definition}\rm}
\newcommand{\edf}{\end{definition}}
\newcommand{\bte}{\begin{theorem}}
\newcommand{\ete}{\end{theorem}}
\newcommand{\bpr}{\begin{proposition}}
\newcommand{\epr}{\end{proposition}}
\newcommand{\ble}{\begin{lemma}}
\newcommand{\ele}{\end{lemma}}
\newcommand{\bco}{\begin{corollary}}
\newcommand{\eco}{\end{corollary}}
\newcommand{\beq}{\begin{equation}}
\newcommand{\eeq}{\end{equation}}
\newcommand{\bdm}{\begin{displaymath}}
\newcommand{\edm}{\end{displaymath}}

\newcommand{\graf}[1]{\left\{\begin{array}{ll}#1\end{array}\right.}

\def\sideremark#1{\ifvmode\leavevmode\fi\vadjust{\vbox to0pt{\vss
 \hbox to 0pt{\hskip\hsize\hskip1em \vbox{\hsize2.1cm\tiny\raggedright\pretolerance10000 \noindent #1\hfill}\hss}\vbox to15pt{\vfil}\vss}}}

\begin{document}
\numberwithin{equation}{section}
\parindent=0pt
\hfuzz=2pt
\frenchspacing

\title[On domains of first kind]{Mean field equations and domains of first kind}

\author[D.B. \& A.M.]{Daniele Bartolucci$^{(1)}$, Andrea Malchiodi$^{(2)}$}


\thanks{$^{(1)}$Daniele Bartolucci, Department of Mathematics, University
of Rome {\it "Tor Vergata"}, \\  Via della ricerca scientifica n.1, 00133 Roma,
Italy. e-mail:bartoluc@mat.uniroma2.it}

\thanks{$^{(2)}$ Andrea Malchiodi, Scuola Normale Superiore, \\ Piazza dei Cavalieri 7, 56126 Pisa, Italy. e-mail:andrea.malchiodi@sns.it}


\maketitle

\setcounter{section}{0}


\setcounter{equation}{0}

\begin{abstract}
	{ In this paper we are interested in understanding the structure of domains of \emph{first and second kind},
		a concept motivated by problems in statistical mechanics. We prove some openness property for
	domains of  first kind with respect to a suitable topology, as well as some sufficient condition for a
simply connected domain to be of first kind in terms of the Fourier coefficients of the Riemann map. Finally, we
show that the set of simply connected domains of first kind is contractible.}
\end{abstract}

\section{Introduction}

We are concerned with the mean field equation,
$$\noindent
\graf{-\Delta u = \lm\dfrac{\dsp e^u}{\ino e^u}\quad \mbox{in}\;\;\om,\\ \\
u=0 \quad \mbox{on}\;\;\pa\om,
}\qquad \qquad \prl
$$
where $\lm\in\R$ and $\om\subset \R^2$ is either any open and bounded domain of class $C^1$ or a bounded simply connected domain,
\underline{regular} according to the following definition { (\cite{CCL})}:

\bdf\label{reg}
{\it Let $\om$ be an open and bounded domain, $\om\subset \mathbb{R}^2$.
We say that $\om$ is \underline{regular}, if its boundary $\partial \om$ is of class $C^2$ but
for a finite number of points $\{Q_1 , . . . , Q_{N_0}\}\subset \partial \om$ such that the following conditions holds at each $Q_j$.\\
(i) The inner angle $\theta_j$ of $\partial \om$ at $Q_j$ satisfies $0 < {\theta_j \neq \pi} < 2\pi$;\\
(ii) At each $Q_j$ there is an univalent conformal map from $B_\delta (Q_j) \cap \overline{\om}$ to
the complex plane $\mathbb{C}$ such that $\partial \om \cap B_\delta (Q_j)$ is mapped to a $C^2$ curve.
 }
\edf

Clearly  any non-degenerate polygon is regular according to this definition.
{ In a slightly different form, based on the concentration/compactenss behavior of minimizers of a mean field variational principle, the following definition was first introduced in \cite{clmp1, clmp2}. Later, a full characterization of the concentration/compactenss behavior of minimizers for simply connected domains has been derived in \cite{CCL}. Finally, the results in
\cite{CCL} have been extended to any connected domain in \cite{BLin3}, thereby establishing the full equivalence with the following:}

\bdf\label{def1} {\it
{ Let $\om\subset\R^2$ be either an open and bounded domain of class $C^1$ or a regular simply connected domain}. We say that $\om$ is of {\bf first kind} if {\rm $\prl$} has no
solution for $\lm=8\pi$. Otherwise $\om$ is said to be of {\bf second kind}. The set of domains of first/second kind will be denoted by $\mathcal{A}_{I}/\mathcal{A}_{II}$ respectively.}
\edf

{ It is worth to point out, as discussed in \cite{clmp2}, that this classification is well understood at least from the
physical point of view.
In the framework of the vortex model of an Euler incompressible flow confined in $\om$,
the Robin function $\gamma_{\ssb}$ is essentially the renormalized free energy of a single vortex.
On domains of first kind the full range of admissible energies $E\in (0,+\ii)$ corresponds to minus the inverse statistical temperature $\lm\in (-\ii,8\pi)$ and as $E\to +\ii$ we have
$\lm\to (8\pi)^-$ and the vorticity of the flow concentrates to a Dirac delta $\delta_{x=q}$, where $q$ is the unique maximum point of the free energy
$\gamma_{\ssb}$. In particular, the equivalence of statistical ensembles holds and the entropy is concave for $E\in (0,+\ii)$.
On domains of second kind the states in the range $\lm\in (-\ii,8\pi)$ describe only a portion of the energy range, say
$E\in (0,E_{8\pi})$, and the peculiar phenomenon of non equivalence of statistical ensembles, see \cite{clmp2}, holds for $E\in (E_{8\pi},+\ii)$
where in particular, as a consequence also of the results in \cite{CCL} and \cite{BLin3}, we have $\lm>8\pi$. Some partial results 
 concerning this problem have been recently obtained in \cite{B2}. Moreover, as $E\to +\ii$, we have $\lm\to (8\pi)^+$ and the vorticity of the flow concentrates to a Dirac delta $\delta_{x=q}$,
where $q$ is a maximum point of the free energy $\gamma_{\ssb}$, see \cite{clmp2}, and a full region of convexity of the entropy is found for $E$
large enough, as first suggested in \cite{B2} and then shown in \cite{bjly2}.}\\

 It is well known that any disk,  say $B_R=B_R(0)$, is of first kind and that, in this particular case, $\mbox{\rm $\prl$}$ admits a solution if and only if $\lm<8\pi$: regular polygons are also of first kind (\cite{CCL}). Symmetric annuli are known to be of second kind, { since a radial solution of $\prl$ exists for any $\lm\in \R$ in this case, see for example \cite{NS}, \cite{BW}}, while $\om=B_R\setminus B_r(x_0)$, with $x_0\neq 0$, is of first kind if $r$ is small enough
 (\cite{BLin3}). Actually this is also an example of a domain of first kind where $\mbox{\rm $\prl$}$ admits solutions also for $\lm>8\pi$. { Indeed, for domains with non-trivial topology, it is well known that for any $N\geq 2$ there are solutions concentrating at $N$ distinct points as
$\lm\to 8\pi N$ \cite{cl1,EGP, KMdP}, as well as solutions for any $\lm\neq 8\pi N$ or for any $\lm$ sufficiently large \cite{cl2, Mal2}.}

It has been proved in \cite{BdM2} that there exists a universal constant $I_c>4\pi$ such that any convex domain whose isoperimetric
ratio $I(\om)$ satisfies $I(\om)>I_c$ is of second kind. Also, if $Q_{a,b}$ is a rectangle whose sides are $1\leq a\leq b<+\ii$ then there exists
$\eta_c\in(0,1)$ such that $Q_{a,b}$ is of second kind if and only if $\frac{a}{b}< \eta_c$, see \cite{CCL}.\\
In particular domains of first kind need not  be symmetric. Let us consider a dumbbell domain $\om_{0,d}$ which is the
union of two disks $B_1, B_2$ of radii $0<r_1\leq r_2$, connected by any smooth thin tube of width $d>0$.
It has been shown in \cite{CCL} that if $r_1<r_2$, then for $d$ small enough
$\om_{0,d}\in \mathcal{A}_{I}$ while if $\om_{0,d}$ is symmetric with respect to the $y$ axis and $r_1=r_2$
then $\om_{0,d}\in \mathcal{A}_{II}$.\\

\brm\label{scale}
{\it By the above discussion, the condition of a domain to be of first or second kind is not conformally invariant.
	However, we recall that $\mbox{\rm $\prl$}$ is scale invariant, that is, $u(x)$ is a solution of $\mbox{\rm $\prl$}$ in $\om$ if and only if
$u(\dt x)$ is a solution of $\mbox{\rm $\prl$}$ in $\frac1\dt \om$. Therefore $\om\in \mathcal{A}_I/\mathcal{A}_{II}$
if and only if $\frac1\dt \om\in \mathcal{A}_I/\mathcal{A}_{II}$ for some $\dt>0$.
}
\erm

We are interested here in a better understanding of the structure of the set of domains of first/second kind. Indeed, besides
the above-mentioned application \cite{clmp2}, this is relevant also for
other problems where, for domains of first kind, one can describe the qualitative behavior of global branches of solutions, see
\cite{BJ,BW} for recent results in this direction.\\\\
Let us recall that, as shown in
\cite{CCL} and \cite{BLin3}, domains of first(second) kind are closed(open) in the $C^1$-domain topology. These results rely on an equivalent
characterization based on the geometric quantity $A_{\ssb}$, see \rife{Dq} and Theorem A below.
Let $G_{\ssb}(x,\,p)$ denote the Green's function of $-\Delta$ with Dirichlet boundary
conditions, uniquely defined by
$$
\left\{
\begin{array}{lll}
-\Delta G_{\ssb}(x,p)&=& \delta_{p}\quad \mbox{in}\quad \om, \\
\hspace{0.3cm}G_{\ssb}(x,p)&=& 0 \quad \mbox{on}\quad \partial \om,
\end{array}\right.
$$
and set
\beq\label{eqn1.2}
    \left\{%
\begin{array}{ll}
    R_{\ssb}(x,\,p)=G_{\ssb}(x,\,p)+\displaystyle\frac{1}{2\pi}\log|x-p|, \\
    \gamma_{\ssb}(p)=R_{\ssb}(p,\,p). \\
\end{array}%
\right.
\eeq
Hence $\gamma_{\ssb}$ denotes the Robin's function relative to $\Omega$ and satisfies,
\beq\label{gmax}
\lim_{p\rightarrow\partial\Omega}\gamma_{\ssb}(p)=-\infty.
\eeq
In view of \rife{gmax}, we see that $\gamma_{\ssb}$ admits at least one critical point, which is its maximum point.
Clearly $q$ is a critical point of $\gamma_{\ssb}$ if and only if $q$ is
a critical point of $R(x,\,q)$ with respect to the $x$ variable,
\beq\label{eqn1.3}
  \left.\nabla_x \gamma_{\ssb}(x)\right|_{x=q}=2\left.\nabla_{x} R_{\ssb}(x,\,q)\right|_{x=q}=0.
\eeq

Let us define,
\beq\label{Dq}
   \pi A_{\ssb}(q)=\lim_{\varepsilon\rightarrow 0}\int_{\Omega\setminus B_\varepsilon (q)}
    \frac{e^{8\pi(R_{\ssb}(x,q)-\gamma_{\ssb}(q))}-1}{|x-q|^4}-\int_{\Omega^{\,c}}\frac{1}{|x-q|^4}.
\eeq

\bigskip

 Here $B_\varepsilon (q)$ denotes the ball of center $q$ and radius
$r$: also, as we will always do later, the standard integration measure has been omitted.

Note that in a neighborhood of $q$,
\beq\label{06.1}
e^{8\pi(R_{\ssb}(x,q)-\gamma_{\ssb}(q))}-1=\sum\limits_{i,j}^{1,2}
a_{ij}(x_i-q_i)(x_j-q_j)+O(|x-q|^3),
\eeq
where, since $R_{\ssb}(x,q)$ is harmonic in $\Omega$,  $a_{11} + a_{22}=0$. In particular, because of
\eqref{06.1}, the limit in \eqref{Dq} always exists and $A_{\ssb}$ is finite.

\medspace

According to some results in \cite{CCL}, \cite{BLin3}
we have the following:

\

{\it \bf Theorem A}
{\it A domain $\om$ of class $C^1$ is of first kind if and only if $\gamma_{\ssb}$ admits a
unique maximum point $q$ and $A_{\ssb}(q)\leq 0$.\\
A simply connected and regular domain $\om$ is of first kind if and only if $\gamma_{\ssb}$ admits a
unique maximum point $q$ and $A_{\ssb}(q)\leq 0$.\\
In particular, in both cases, if $\gamma_{\ssb}$ admits a critical point $q$ such that $A_{\ssb}(q)\leq 0$,
then $q$ is the unique maximum point, it
is a non-degenerate critical point of $\gamma_{\ssb}$ and $\om$ is of first kind.}\\

Although  not stated in this form in the cited references, still Theorem A is a trivial consequence of the results obtained therein. We remark that
the proof of  Theorem A crucially relies on the uniqueness and non-degeneracy of solutions of $\prl$ for $\lm=8\pi$, see \cite{CCL} and \cite{BLin3}.
If we miss the regularity assumptions on $\om$ in the claim, then we do not know much about this point and in particular about the validity of
Theorem A.
\bigskip

Next we focus on the case where $\om$ is simply connected and use complex notation. For a fixed $q\in\om$, let us denote by
$\mathbb{D}=\{z\in \ci\,:\,|z|<1\}$ and by $g_q:\om \mapsto \mathbb{D}$ the Riemann map
satisfying $g_q(q)=0$, $g_q^{'}(q)>0$.
Let $f_q:\mathbb{D}\mapsto \om$ be the inverse map, which satisfies $f_q(0)=q$ and
$g_q=\fq^{-1}:\om \mapsto \mathbb{D}$. Next,  setting $w=f_q(z)$, we find that
$$
R_{\ssb}(w,q)=G_{\ssb}(w,q)+\frac{1}{2\pi}\log{|w-q\,|}=-\frac{1}{2\pi}\log{\frac{|g_q(w)|}{|w-q\,|}}=
\frac{1}{2\pi}\log{\frac{|\fq(z)-\fq(0)\,|}{|z|}},
$$
and it is well known that the Robin function takes the form
\beq\label{gammaex}
\gamma_{\ssb}(w)=R_{\ssb}(w,w)=\frac{1}{2\pi}\log{\frac{(1-|g_q(q)|^2)}{|g_q^{'}(w)|}}=
\frac{1}{2\pi}\log{(1-|z|^2)|\fq^{'}(z)|}.
\eeq
Next, let us consider the power series relative to $\fq$,
\beq\label{fq1}
\fq(z)=q+a_1z+\sum\limits_{n=2}^{+\ii}a_nz^n,\quad |z|<1,
\eeq
where we used that $a_1=\fq^{'}(0)>0$, since by assumption $g_q^{'}(q)>0$. Therefore we see that
$$
\gamma_{\ssb}(q)=\frac{1}{2\pi}\log(|a_1|),
$$
and
\beq\label{a2}
\left.\frac{\pa}{\pa z}R_{\ssb}(z,q)\right|_{z=q}=\frac{a_2}{4\pi a_1^2}=\frac{a_2}{4\pi},
\eeq
whence \eqref{eqn1.3} is equivalent to $a_2=0$. At this point one can prove (see \cite{CCL}) that,

\beq\label{ADq}
D_{\ssb}(q):=|a_1|A_{\ssb}(q)=-|a_1|^2+\sum\limits_{n=3}^{+\ii}\frac{n^2}{n-2}|a_n|^2,
\eeq
which is well defined and convergent since $a_2=0$ and $|\om|=\pi\sum\limits_{n=1}^{+\ii}n |a_n|^2$ by the Area theorem.
We will use when needed the fact that $D_{\ssb}$ and $A_{\ssb}$ share the same sign without further comments.
In particular, any bounded and simply connected domain admits a Riemann map whose series expansion takes the form,
\beq\label{fq2}
\fq(z)=q+a_1z+\sum\limits_{n=3}^{+\ii}a_nz^n,\quad |z|< 1,\qquad a_1>0,
\eeq
where $q$ is a critical point of $\gamma_{\ssb}$.\\
For a pair $\om,\om_0$ of class $C^1$($C^{0,1}$), we will denote by $d_1(\om,\om_0)(d_{0,1}(\om,\om_0))$
the distance in  the $C^1$($C^{0,1}$)-domain topology, see Section \ref{secP} for more details.

\medspace

We then have the following result:

\

\bte\label{thm2}$\,$
$(i)$ Let $\om$ be a domain of first kind of class $C^1$ with $A_{\ssb}(q)<0$.
Then there exists $\eps_0>0$ such that if $\om_0$ satisfies $d_1(\om,\om_0)<\eps_0$, then
$\om_0\in \mathcal{A}_{I}$ and $A_{\sscp \om_0}(q_0)<0$.\\
$(ii)$ Let $\om$ be a simply connected and regular domain of first kind with $A_{\ssb}(q)<0$.
Then there exists $\eps_0>0$ such that if $\om_0$ is a simply connected and regular domain and
$d_{0,1}(\om,\om_0)<\eps_0$, then $\om_0\in \mathcal{A}_{I}$ and $A_{\sscp \om_0}(q_0)<0$.
\ete

\bigskip
It is interesting that $\mathcal{A}_I\cap \{A_{\ssb}(q)<0\}$ is open in the $C^1$-topology.
Obviously the problem is more subtle for the $C^{0,1}$-domain topology. Indeed our proof of Theorem \ref{thm2}
crucially relies on Theorem {A}, that is, on a characterization of domains of first kind, which unfortunately is not known for a general domain of class $C^{0,1}$.\\\\

Our next result is a coefficient-based sufficient condition that defines an open region of starlike domains
in $\mathcal{A}_I$ which contains all disks.
Let $\om$ be a simply connected domain,
$q$ a critical point of $\gamma_{\ssb}$ and $\fq$ the Riemann map defined in \rife{fq2}.
Let us denote by $\mathcal{S}_{I}$ the subset of those $\om$ such that $\om=\fq(\mb{D})$ with $\fq$ as in \rife{fq2} and
\beq\label{nan}
\sum\limits_{n=3}^{+\ii}n|a_n|< |a_1|,
\eeq
and by $\partial \mathcal{S}_{I}$ the subset of those $\om$ such that,
\beq\label{nanbdy}
\sum\limits_{n=3}^{+\ii}n|a_n|= |a_1|.
\eeq
Then we have the following:
\bte\label{thm3}
If $\om\in \mathcal{S}_{I}$, then $\om$ is starlike,  of class $C^1$ and of first kind with $D_{\ssb}(q)<0$.
In particular, if $\om\in \mathcal{S}_{I}$
and $\{a_n\}_{n\in\N}$ are the coefficients of {\rm \rife{fq2}}, then for any
continuous (w.r.t. the $\ell_\ii$-topology) map $\un{a}(t)=(a_1,0,a_3(t),a_4(t),\cdots)$, $t\in[0,1]$, satisfying $a_n(0)=0$,  $a_n(1)=a_n$ and
$|a_n(t)|\leq |a_n|$  for any $n\geq 3$, then
$$
f(z,t)=q+a_1z+\sum\limits_{n=3}^{+\ii}a_n(t)z^n,\quad |z|\leq 1, \quad t\in [0,1],
$$
is a jointly continuous family of univalent and starlike maps $f(z,t)$, $z\in \ov{\mathbb{D}}$, $t\in [0,1]$ such that $\ov{\om_t}=f(\ov{\mathbb{D}},t)$
satisfies $\om_t\in \mathcal{S}_{I}$ for any $t$, $\om_0=q+a_1\mathbb{D}$, $\om_1=\om$. Therefore, in particular $\om_t$ is of
first kind and $D_{\sscp \om_t}(q)<0$ for any $t\in[0,1]$.\\
Moreover, the same conclusion holds for $\om_t$ with $t\in[0,1)$,
whenever $\om\in \pa\mathcal{S}_{I}$ is regular and if at least one $|a_n(t)|$ is strictly
increasing in a left neighborhood of $t=1$.
\ete

Please observe that, since \rife{nan} implies that $\om$ is of class $C^1$,  it also follows from Theorem \ref{thm2}
that $\mathcal{S}_{I}$ is open. However the interest of Theorem \ref{thm3} relies in the fact that, as already mentioned above,
 it is not true that
any starlike or either convex domain is of first kind.\\
Actually, since by the Area formula? the formal series built with $a_n(t)$ as in Theorem \ref{thm3}, that is
$D_{\ssb}(q)=-|a_1(t)|^2+\sum\limits_{n=3}^{+\ii}\frac{n^2}{n-2}|a_n(t)|^2$, is  convergent, then
it is tempting to try to adopt the same argument to prove the path-connectedness of the full set
of $C^{1}$ or regular domains of first kind.\\ Unfortunately this argument fails in general.
Indeed, on one side if we miss \rife{nan} then it is not anymore guaranteed that $f(\mathbb{D})$ is either $C^{1}$ or even just
regular in the sense of Definition \ref{reg}. For example, as observed in \cite{CCL}, $f_3(z)=z+\frac13 z^3$ is univalent in $\mb{D}$
and satisfies $D_{\ssb}(0)=\sum\limits_{n=3}^{+\ii}\frac{n^2}{n-2}|a_n|^2-|a_1|^2= 9|a_3|^2-|a_1|^2=0$,
$\sum\limits_{n=3}^{+\ii}n|a_n|=3|a_3|=1=|a_1|$, but at the same time $f_3^{'}(\pm i)=0$ and the domain $f_3(\mathbb{D})\in \pa \mathcal{S}_{I}$ has cusps
(inner angle $2\pi$) at its boundary points $\pm \frac23 i\in \pa\om$. In particular $f_3(\mathbb{D})$ is not regular and then
we cannot apply Theorem A above, whence the full argument breaks down.
Actually we see in this way that \rife{nan} is sharp as far as we are concerned with the regularity of the domain. However \rife{nan}
is not necessary for a domain to be of first kind, as we illustrate with an explicit example in Appendix III.\\
On the other
side it is neither true that if $D_{\ssb}(q)=-|a_1(t)|^2+\sum\limits_{n=3}^{+\ii}\frac{n^2}{n-2}|a_n(t)|^2$ is negative then
$f$ as in \rife{fq2} is univalent. For example the sequence
$\un{a}=(1,0,\frac{t}{3},\frac{t}{4},0,0,\cdots)$ corresponds to the holomorphic function $f(z)=z+\frac{t}{3}z^3+\frac{t}{4}z^4$,
which for $t\in [\frac{7}{10},\frac{4}{5}]$ is readily seen to be not univalent in $\mathbb{D}$, although $D_{\ssb}(0)$ is
 convergent and strictly negative.
\\\\

However we can show that a particularly simple choice of the $a_n(t)$'s in Theorem \ref{thm3} does the job for any simply connected
domain of first kind. We recall that a domain $\om$ is said to be \emph{analytic} if $\pa\om=f(\pa \mb{D})$ for some $f$ univalent in a full
open neighborhood of $\mb{D}$. Then we have,
\bte\label{pathc} Let $\om$ be a simply connected domain of first kind, either regular or of class $C^{1}$, and let $\fq$ be
the Riemann map normalized as in {\rm \rife{fq2}}.
Then
\beq\label{fqt}
f(z,t)=tq+\frac{\fq(tz)-q}{t},\quad |z|\leq 1, \quad t\in [-1,1],
\eeq
is jointly continuous in $\ov{\mathbb{D}}\times [-1,1]$, jointly analytic in
$(z,t)\in \mathbb{D}\times t\in(-1,1)$, $\ov{\om_t}=f(\ov{\mathbb{D}},t)$ is an analytic domain for any $t\in [0,1)$
and satisfies $\om_t\in \mathcal{A}_{I}$ for any $t\in[0,1]$,  $\om_1=\om$ and $\om_0=\mathbb{D}$.\\
In particular, the set of simply connected $C^1$ domains of first kind is contractible, while the set of simply connected
regular domains of first kind is simply connected w.r.t. the $C^{0,1}$-topology.
\ete

What the proof shows is that $f(z,t)$ defines a deformation retract of the identity in the subspace of $C^1$ domains of first kind.
As remarked right after Theorem \ref{thm2}, the situation for regular domains is more delicate, which is why we come up with a
weaker result in this case. In conclusion, as a consequence of Theorems \ref{thm2} and \ref{pathc}, we have the following,
\bco\label{co7} The set of simply connected domains of first kind of class $C^1$ is a contractible set with non empty interior with respect
to the $C^1$-topology.\\
The set of regular and simply connected domains of first kind is a simply connected set with respect
to the $C^{0,1}$-topology.
\eco
\bigskip

It is an interesting open problem to understand how the topology of $\mathcal{A}_{I}$ is
affected by the topology of the underlying domains.

\bigskip
\bigskip

This paper is organized as follows. In Section \ref{secP} we define the distances and topologies used
in the introduction and list some known results which will be needed in the sequel. In Section \ref{secthm2} we prove Theorem \ref{thm2}.
In Section \ref{secthm3} we prove Theorems \ref{thm3} and \ref{pathc}. Some technical results and an example are discussed
in the Appendices.\\

\medspace

\begin{center}
	{\bf Acknowledgements}
\end{center}

A.M. is supported by the project {\em Geometric problems with loss of compactness} from Scuola Normale Superiore and
by MIUR Bando PRIN 2015 2015KB9WPT$_{001}$.  He is also member of GNAMPA as part of INdAM.\\
D.B. is partially supported by:
MIUR Bando PRIN 2015 2015KB9WPT$_{001}$, Consolidate the Foundations project 2015 (sponsored by Univ. of Rome "Tor Vergata") "{\em Nonlinear Differential Problems and their Applications}",
S.E.E.A. project 2018 (sponsored by Univ. of Rome "Tor Vergata"),
MIUR Excellence Department Project awarded to the Department of Mathematics, Univ. of Rome Tor Vergata, CUP E83C18000100006.

\bigskip
\bigskip

\section{Preliminaries}\label{secP}

Let us now introduce some useful definitions and distances between domains.

\bdf\label{C1} {\it
A domain $\om$ is of class $C^{k}(C^{0,1})$, $k\geq2$, if for each $x_0\in \pa\om$ there exists a ball $B=B_r(x_0)$
and a one-to-one map $\Phi: B\mapsto U\subset \R^2$ such that $\Phi\in C^{k}(B)(C^{0,1}(B)),
\Phi^{-1}\in C^{k}(U)(C^{0,1}(U))$ and the following holds:
$$
\Phi(\om\cap B)\subset \R^2_+\quad \mbox{ and }\quad \Phi(\om\cap B)\subset \pa\R^2_+.
$$

It is well known {\rm(}see for example \cite{H}{\rm)} that this is equivalent to the existence of $r>0$ and $M>0$ such that, given any ball
$B_r(x_0)$, $x_0\in \R^2$ then, after suitable rotation and translations, it holds:
$$
\om\cap B=\{(x_1,x_2)\,:\,x_2<\phi(x_1)\}\cap B\quad \mbox{ and } \quad \pa\om\cap B=\{(x_1,x_2)\,:\,x_2=\phi(x_1)\}\cap B,
$$

for some $\phi\in C^{k}(\R)(C^{0,1}(\R))$.}
\edf

We will also need some classical results about extensions up to the boundary of Riemann maps.

\brm\label{rempom} {\it If $f:\mb{D}\to \om$ is univalent and if $\pa\om$ is the support of any
rectifiable Jordan curve, then $f$ admits a continuous and univalent extension on $\ov{\mb{D}}$, see for example Theorem 9.1 in \cite{pom2}.
Moreover, if $\om$ is of class $C^{1}$, then, in view of Definition \ref{C1}, it is not difficult to see that
$\pa\om$ admits a $C^1$-parametrization $w(t)$, $t\in [0,2\pi]$, such that $w^{'}(t)\neq 0$, $t\in [0,2\pi]$. As a consequence, see Theorem 3.5 in
\cite{pom}, $f^{'}$ admits a continuous extension on $\ov{\mb{D}}$, with $f^{'}(z)\neq 0$ in $\ov{\mb{D}}$.  }
\erm

Next we will use the following definition of distance in the set $\mathcal{M}_1(\om)$ of bounded $C^1$ domains which are
$C^{1}$-diffeomorphic to a
given bounded domain $\om$. This is a particular case of a more general definition first introduced in \cite{Mich}, see also \cite{H}.
Let $\om_1$, $\om_2\in\mathcal{M}_1(\om)$: then we define

\beq\label{d1}
d_1(\om_1,\om_2)=\inf\left\{\sum\limits_{k=1}^N(\|h_k-I\|_{C^1}+\|h^{-1}_k-I\|_{C^1})\,|\,h_1\circ h_2 \circ \cdots \circ h_N(\om_1)
=\om_2\right\},
\eeq
where $I:\R^2\to \R^2$ is the identity and the infimum is taken over $N \in \N$ and all $C^1$-diffeomorphisms $h_k:\R^2\to \R^2$ such that
$D^{j}(h_k(x)-x)\to 0$ as $|x|\to +\ii$, $j=0,1$.\\
Equipped with this metric, $\mathcal{M}_1(\om)$ is complete and separable. A neighborhood of $\om_1\in\mathcal{M}_1(\om)$
in the induced topology contains a neighborhood of the form $$
\{H(\om_1)\,:\,H\in C^{1}(\R^2;\R^2),\,\|H-I\|_{C^1}<\eps,\,D^{j}(H(x)-x)\to 0,\,|x|\to +\ii,\,j=0,1\},
$$
which in turn contains a ball $\{\om_2\,:\,d_1(\om_2,\om_1)<\delta\}$, for some $\delta>0$, see Appendix A.2 in \cite{H} for
proofs. This is the distance $d_1$ and the $C^1$-domain topology which we refer to in Theorems \ref{thm2},
\ref{pathc} and Corollary \ref{co7}. In particular the set of all open and bounded domains of class $C^1$ splits into equivalence
classes with respect to the relation $\om_1\backsim\om_2 \Leftrightarrow \om_2=h(\om_1)$ for some $C^1$-diffeomorphism $h$. Obviously
one of these equivalence classes is the subset of simply connected domains of class $C^1$.\\
Concerning simply connected and
regular domains (see Definition \ref{reg}), we first observe that, in view of Remark \ref{rempom},
any such a domain $\om$ can be mapped one-to-one onto $\mb{D}$ via a Riemann map $f:\mb{D}\to \om$,
which admits a one-to-one and continuous extension on $\ov{\mb{D}}$. Let $Q_j\in \pa{\om}$ be any corner
and assume without loss of generality that $Q_j=0=f(1)$. Then,
by Theorem {3.9} in \cite{pom}, either $\theta_j=\pi\alpha\in (0,\pi)$, and then $|f^{'}(z)|\leq C |z-1|^{1-\al}$ for
$z\in \ov{\mb{D}}\cap B_r(1)$, or $\theta_j=\pi\alpha\in (\pi,2\pi)$, and then $\frac{|f(z)|}{|z-1|}\leq C |z-1|^{\al-1}$ for
$z\in \ov{\mb{D}}\cap B_r(1)$, for suitable $r>0$. Therefore, it is not difficult to see that any regular domain is in particular
of class $C^{0,1}$ in the sense defined above. As a consequence the distance of any two regular domains, which we will denote
by $d_{0,1}(\om_1,\om_2)$,
can be defined just by replacing $C^{1}$ with
$C^{0,1}$ in \rife{d1}. A generalization of the notion introduced in \cite{Mich}, with a proof of the fact that $d_{0,1}$ is indeed a
well defined metric, can be found in Chapter 3 of \cite{DZ}.

\bigskip

\section{Proof of Theorem \ref{thm2}}\label{secthm2}
In this Section we prove Theorem \ref{thm2},  first discussing $(ii)$. The proof of $(i)$ when $\om$ is simply connected follows exactly by the same argument but it is easier
and we omit it here to avoid repetitions. Then we will be back to $(i)$ for general connected but not simply connected
domains $\om$ of class $C^1$.\\
We can assume without loss of generality that $q=0$ and in particular,
by Remark \ref{scale}, that $a_{1}=\fq^{'}(0)=1$. By Remark \ref{rempom} the Riemann map \rife{fq2} can be extended to a
continuous and univalent map
$\fq(\ov{\mathbb{D}})=\ov{\om}$ which then takes the form,
$$
\fq(z)=z+\sum\limits_{n=3}^{+\ii}a_n z^n,\quad z\in \ov{\mathbb{D}}.
$$
We argue by contradiction. If the claim were false then we could find a sequence of regular domains $\om_k$ such that $d_{0,1}(\om,\om_k)\to 0$ in the $C^{0,1}$-topology and
$\om_k\in \mathcal{A}_{II}$ for any $k$. Let
$\gamma_{\sscp \om_k}$ denote the Robin function of $\om_k$.
Since $d_{0,1}(\om,\om_k)\to 0$ and each $\om_k$ is regular, then it is not difficult to see that $\om_k\to \om$ in the sense of
Caratheodory's kernel convergence \cite{pom2}. Therefore, for any $k$ large enough, since $q=0$ is an interior point of $\om$, then
we have $q=0\in \om_k$ and then
we can define $\widehat{f}_{\sscp k}:\ov{\mathbb{D}}\to\ov{\om}_k$ to be a sequence
of univalent maps which satisfy
$\widehat{f}_{\sscp k}(\ov{\mathbb{D}})=\ov{\om}_k$ and hence
$$
\widehat{f}_{k}(z)=\widehat{a_{1,k}}z+\sum\limits_{n=2}^{+\ii}\widehat{a_{n,k}}z^n\quad z\in \ov{\mathbb{D}},\qquad
\widehat{a_{1,k}}>0,\,\forall\,k\in \N.
$$
By the Caratheodory kernel Theorem (\cite{pom2}, Theorem 1.8) we conclude that
$\widehat{f}_{k}\to \fq$ locally uniformly and then also in $C_{\rm loc}^3(\mathbb{D})$.\\
Next observe that, by Theorem {A}, $q=0=\fq(0)$ is the unique and non-degenerate maximum point of $\gamma_{\ssb}$.
As a consequence of \rife{gammaex}, we see that $\gamma_{\sscp \om_k}\to \gamma_{\ssb}$ in $C_{\rm loc}^2(D)$ and then,
for $k$ large enough, $\gamma_{\sscp \om_k}$ has a unique and non-degenerate maximum point, which we denote by $q_k$ and satisfies
$q_k\to q=0$. At this point we can define
\beq\label{gk}
{f}_{\sscp k,q_k}(z)=q_k+a_{1,k}z+\sum\limits_{n=2}^{+\ii}{a_{n,k}}z^n\quad z\in \ov{\mathbb{D}},\qquad {a_{1,k}}>0,\,\forall\,k,
\eeq
to be the sequence of univalent functions which satisfies ${f}_{\sscp k,q_k}(\ov{\mathbb{D}})=\ov{\om}_{k}$, ${f}_{\sscp k,q_k}(0)=q_k$ and
$f^{'}_{\sscp k,q_k}(0)>0$.  As a consequence of \rife{a2} we find that $a_{2,k}=0$. In particular, by using once more the kernel Theorem,
we find that ${f}_{\sscp k,q_k}\to \fq$ locally uniformly in $\mathbb{D}$.\\
Next, let us recall that a sequence of compact domains $U_k$ is said to be \emph{uniformly locally connected} if  for every $\eps>0$ there exists
$\dt>0$ such that if $a_k,b_k \in \om_k$ and $|a_k-b_k|<\dt$ there exists connected compact sets $B_k\subseteq U_k$ such that
$a_k,b_k\in B_k$ and diam$(B_k)<\eps$. Since $d_{0,1}(\om, {\om}_k)\to 0$ and each ${\om}_k$ is regular, then it can be shown that:\\\\
{\bf Claim:}
$\ov{\om}_k$ is uniformly locally connected.\\
See Appendix I for a proof of this fact.\\\\
Therefore, by Theorem 9.11 in \cite{pom2}, we conclude that ${f}_{\sscp k,q_k}\to \fq$ uniformly in $\ov{\mathbb{D}}$.
Thus, since $\|{f}_{\sscp k,q_k}-\fq\|_{\ii}\to 0$, then by  Cauchy's representation formula we find that
$$
\sup\limits_{n\in\N} |a_{n,k}- a_n|\leq \|{f}_{\sscp k,q_k}-\fq\|_{\ii}\to 0, \mbox{ as }k\to +\ii.
$$

At this point, since $D_{\ssb}(q)<0$, we have
$$
\lim\limits_{N\to +\ii}\lim\limits_{k\to +\ii} \sum\limits_{n=3}^{N}\frac{n^2}{n-2}|{a_{n,k}}|^2=
\sum\limits_{n=3}^{+\ii}\frac{n^2}{n-2}|{a_{n}}|^2=|a_1|^2-\sg=1-\sg,
$$
for some $\sg\in (0,1)$, while on the other side, by the uniform convergence of $a_{n,k}$, we find that
$$
1-\sg= \lim\limits_{k\to +\ii}\lim\limits_{N\to +\ii} \sum\limits_{n=3}^{N}\frac{n^2}{n-2}|{a_{n,k}}|^2=
\lim\limits_{k\to +\ii} \sum\limits_{n=3}^{+\ii}\frac{n^2}{n-2}|{a_{n,k}}|^2\geq \lim\limits_{k\to +\ii}|a_{1,k}|^2=1,
$$
where the last inequality follows from Theorem {A}, which implies that, since $\om_{k}\in\mathcal{A}_{II}$  for any $k$, then
$\sum\limits_{n=3}^{+\ii}\frac{n^2}{n-2}|{a_{n,k}}|^2> |a_{1,k}|^2$ for any $k$.
This contradiction shows that there exists $\eps_0>0$ such that if $\om_0$ is regular and $d_{0,1}(\om,\om_0)<\eps_0$ then
$\om_0\in \mathcal{A}_{I}$.\\
Taking a smaller $\eps_0$ if necessary, the same argument shows that  $D_{\sscp \om_0}(q_0)<0$ as well and we
skip this part of the proof to avoid repetitions.

\bigskip

We are left with the proof of $(i)$ in the case when $\om$ is not simply connected.\\
We will denote by $C>0$ a uniform positive constant whose value may change from line to line.
By Theorem A, $\om$ has a unique and non-degenerate maximum point $q$.
If the claim were false then we could find a sequence of $C^{1}$ domains $\om_k$ such that $d_{1}(\om,\om_k)\to 0$
and $\om_k\in \mathcal{A}_{II}$ for any $k$. Let $R_k(x,y)$ be the regular part of the Green function for $\om_k$ as defined
in \rife{eqn1.2} and $\gamma_k$ its Robin function. It can be shown that, see Appendix II,  for any fixed
$y\in \om$ we have
\beq\label{A.a1}
R_{k}(x,y)\to R_{\ssb}(x,y) \mbox{ in } C^{3}_{\rm loc}({\om})
\eeq
and
\beq\label{A.b1}
\gamma_k\to \gamma_{\ssb} \mbox{ in } C^{2}_{\rm loc}(\om).
\eeq

Therefore, for $k$ large, $\gamma_k$ will have a unique and non-degenerate maximum point $q_k\to q$. We can assume for the moment without loss of generality
that $q_k=0$ for any $k$. By assumption, for fixed $k$, we have,
$$
    A_{k}(0):=\lim_{\varepsilon\rightarrow 0}\int_{\om_k\setminus B_\varepsilon (0)}
    \frac{e^{8\pi(R_{k}(x,0)-\gamma_{k}(0))}-1}{|x|^4}-\int_{\om^{\,c}}\frac{1}{|x|^4}>0,
$$
and then in particular
$$
\liminf\limits_{k\to+\ii} A_{k}(0)\geq 0.
$$
We will obtain a contradiction by showing that there exists $\sg>0$ small enough such that,
\beq\label{contr}
\limsup\limits_{k\to+\ii}  A_{k}(0)\leq -\sg.
\eeq
Clearly there exists $d>0$ small enough such that $B_d(0)\subset \subset \om_k$ for any $k$ large.
By \rife{06.1} we have
$$
e^{8\pi(R_{k}(x,0)-\gamma_{k}(0))}-1=\sum\limits_{i,j}^{1,2}
a_{k,ij}x_ix_j+c_k(x),\;x\in \ov{B_d(0)},
$$
where
$$
a_{k,ij}=\pa_{x_i,x_j}R_k(x,0)\to \pa_{x_i,x_j}R(x,0)=a_{ij},\quad i,j=1,2,
$$
and by \eqref{A.a1} the reminders $c_k$ satisfy
$$
|c_k(x)|\leq C|x|^3,\; 
\quad \mbox{\rm for $k$ large}.
$$
As a consequence, for any $\varepsilon \leq \frac{d}{4}$, we find that,
$$
\left|\int_{B_d(0)\setminus B_\varepsilon (0)}
    \frac{e^{8\pi(R_{k}(x,0)-\gamma_{k}(0))}-1}{|x|^4}\right|\leq
    \left|\int_{B_d(0)\setminus B_\varepsilon (0)}\frac{\sum a_{k,ij}x_ix_j+c_k(x)}{|x|^4}\right|=
$$
\beq\label{hat0}
\left|\int_{B_d(0)\setminus B_\varepsilon (0)}\frac{2a_{k,12}x_1x_2+c_k(x)}{|x|^4}\right|=
\left|\int_{B_d(0)\setminus B_\varepsilon (0)}\frac{ c_k(x)}{|x|^4}\right|\leq Cd ,\; 
\quad \mbox{\rm for $k$ large},
\eeq
where we used the symmetry of the domain and the fact that, since $R_{k}(x,0)$ is harmonic in $\om$,  $a_{k,11} + a_{k,22}=0$ for any $k$.
To simplify the evaluation let us set
$$
h_k(x)=\frac{e^{8\pi(R_{k}(x,0)-\gamma_{k}(0))}-1}{|x|^4},\,x\in\om_k,\quad
h(x)=\frac{e^{8\pi(R_{\ssb}(x,0)-\gamma_{\ssb}(0))}-1}{|x|^4},\,x\in\om.
$$

Since $R_{k}(x,q_k)\to R_{\ssb}(x,q)$ locally uniformly in $\om$, and $\gamma_k(q_k)\to \gamma_k(q)$,
then, for any open and relatively compact subset $\widehat{\om}\subset\subset \om$, we have:
\beq\label{hat2}
\widehat{\om}\subset\subset \om_k,\,\mbox{ for $k$ large and }\int_{\widehat{\om}\setminus B_d(q_k)}
    h_k(x)\to \int_{\widehat{\om}\setminus B_d(q)} h(x),\;k\to +\ii.
\eeq

Since the symmetric difference $\om_k \Delta \om\to \emptyset$ as $k\to +\ii$, then for any $\dt>0$ we can choose
an open and relatively compact subset $\widehat{\om}_{\dt}$ as in \rife{hat2} which also satisfies,
$$
\om_k \setminus \widehat{\om}_{\dt} \subset B_{2\dt}(\pa \om_k)\mbox{ and } \om \setminus \widehat{\om}_{\dt} \subset B_{2\dt}(\pa \om)\;\mbox{\rm for $k$ large},
$$
and then  in particular,
\beq\label{hat1.2}
\int_{(\om_k\setminus B_d(q_k))\setminus (\widehat{\om}_{\dt}\setminus B_d(q_k))}
    \left|h_k(x)\right|\leq \int_{B_{2\dt}(\pa \om_k)}
    \left|h_k(x)\right|\leq C |\pa \om_k|\dt, \;\mbox{\rm for $k$ large},
\eeq
\beq\label{hat1.3}
\int_{(\om\setminus B_d(q))\setminus (\widehat{\om}_{\dt}\setminus B_d(q))}
    \left|h(x)\right|\leq \int_{B_{2\dt}(\pa \om)}
    \left|h(x)\right|\leq C |\pa \om|\dt, \;\mbox{\rm for $k$ large},
\eeq
where we used the uniform bound
\beq\label{A.c1}
|h_k(x)|\leq C, \; x\in B_{2\dt}(\pa \om_k)\cap \om_k,
\eeq
see Appendix II. Thus we can estimate,
$$
\int_{\om_k\setminus B_{\varepsilon}(q_k)}h_k(x)-\int_{\om\setminus B_{\varepsilon}(q)}h(x)\leq
\int_{B_d(q_k)\setminus B_{\varepsilon}(q_k)}|h_k(x)|+\int_{B_d(q)\setminus B_{\varepsilon}(q)}|h(x)|+
$$
$$
\int_{(\om_k\setminus B_d(q_k))\setminus (\widehat{\om}_{\dt}\setminus B_d(q_k))}|h_k(x)|+
\int_{(\om\setminus B_d(q))\setminus (\widehat{\om}_{\dt}\setminus B_d(q))}|h(x)|+
\left|\int_{\widehat{\om}_{\dt}\setminus B_d(q_k)}h_k(x)-\int_{\widehat{\om}_{\dt}\setminus B_d(q)}h(x)\right|\leq
$$
\beq\label{hat4}
Cd+C\dt+\left|\int_{\widehat{\om}_{\dt}\setminus B_d(q_k)}h_k(x)-\int_{\widehat{\om}_{\dt}\setminus B_d(q)}h(x)\right|,
\eeq
where we used \rife{hat0},\rife{hat1.2},\rife{hat1.3}.
At this point let us fix $\sg>0$ such that $A_{\ssb}(q)=-5\sg$ and then choose $\dt$ and $d$ such that
$Cd+C\dt<\sg$. For any $\dt$ and $d$  fixed in this way, by \rife{hat2} we have,
$$
\left|\int_{\widehat{\om}_{\dt}\setminus B_d(q_k)}h_k(x)-\int_{\widehat{\om}_{\dt}\setminus B_d(q)}h(x)\right|< \sg,
$$
for any $k$ large enough. In particular, for $k$ large we also have,
$$
\left|\int_{\om_k^c}\frac{1}{|x-q_k|^4}-\int_{\om^c}\frac{1}{|x-q|^4}\right|<\sg.
$$

Finally, we can choose $\varepsilon_0>0$ small enough to guarantee that
$$
\int_{\om\setminus B_{\varepsilon}(q)}h(x)-\int_{\om^c}\frac{1}{|x-q|^4}\leq -4\sg,
$$
for any $\varepsilon< \varepsilon_0$. Plugging these estimates together with \rife{hat4} we conclude that
$$
\int_{\om_k\setminus B_{\varepsilon}(q_k)}h_k(x)-\int_{\om_k^c}\frac{1}{|x-q_k|^4}\leq
\int_{\om\setminus B_{\varepsilon}(q)}h(x)-\int_{\om^c}\frac{1}{|x-q|^4}+3\sg\leq -\sg,
$$
for any $k$ large enough and for any $\varepsilon< \varepsilon_0$. As a consequence we conclude that \rife{contr} holds,
which is the desired contradiction.

\bigskip

\section{Proofs of Theorems \ref{thm3} and \ref{pathc}}\label{secthm3}
In this Section we prove Theorems \ref{thm3} and \ref{pathc}.\\

{\it  Proof of Theorem \ref{thm3}.}
By Remark \ref{scale} we can assume without loss of generality that $a_1=\fq^{'}(0)=1$ and  after  a translation we can also assume that
$q=0$. Therefore \rife{fq2} takes the form 
$$
\fq(z)=z+\sum\limits_{n=3}^{+\ii}a_nz^n,\quad |z|< 1.
$$

It is well known (see for example \cite{pom2} p.44) that if $\sum\limits_{n=2}^{+\ii}n |a_n|\leq 1$ then $\fq^{'}$ has positive real part
and $\fq$ is univalent and starlike in $\mathbb{D}$. Therefore, since $a_2=0$, by \rife{nan} $\fq$ is infact univalent and starlike.
Setting $h(\theta)=\fq(e^{i\theta})$, $\theta \in  [0,2\pi]$, we have
$$
h(\theta)= e^{i\theta}+\sum\limits_{n=3}^{+\ii}a_ne^{in\theta},
$$
and once more \rife{nan} shows that the real and imaginary parts of $h$ have continuous first derivative in $[0,2\pi]$
satisfying $\|h^{'}\|_{\ii}\leq 2$.
On the other side we also have,
$$
|h^{'}(\theta)|\geq 1-\sum\limits_{n=3}^{+\ii}n|a_n|>0,
$$
once more by \rife{nan}. Therefore $h(\theta)$ is a $C^1$ curve and then in particular $\om$ is of class $C^1$. Since
$\sum\limits_{n=3}^{+\ii}n |a_n|< 1$, then $\sum\limits_{n=3}^{+\ii}\frac{n^2}{n-2} |a_n|^2\leq
\sum\limits_{n=3}^{+\ii}n^2 |a_n|^2<  1$. Therefore $D_{\ssb} (q)< 0$ and hence $\om\in \mathcal{A}_I$ by Theorem A.\\
In particular we have shown that any $\om\in \mathcal{S}_{I}$ is a $C^1$ domain of first kind with $D_{\ssb} (q)< 0$.
At this point we define,
$$
f(z;t)=z+\sum\limits_{n=3}^{+\ii}a_n(t)z^n,\quad |z|<1, \quad t\in [0,1],
$$
where each $a_n(t)$ is continuous in $[0,1]$, $a_n(0)=0$,  $a_n(1)=a_n$ and
$|a_n(t)|\leq |a_n|$  for any $n\geq 3$. It is easy at this point to see that $\om_t=f({\mathbb{D}},t)\in \mathcal{S}_{I}$ for any $t$.
In particular, by Remark \ref{rempom} $f$ admits a continuous and univalent extension
$f(\ov{\mathbb{D}},1)=\ov{\om}$ and $f(\ov{\mathbb{D}},0)=q+\ov{\mathbb{D}}$, as claimed.\\
Finally, if $\om\in \pa \mathcal{S}_I$, then the same argument shows that $\om_t=f(\mathbb{D},t)\in \mathcal{S}_I$ for any $t\in [0,1)$,
whenever at least one $|a_n(t)|$ is strictly increasing for $t\simeq 1^-$, and the conclusion follows in this case as well. \finedim

\bigskip
\bigskip

{\it  Proof of Theorem \ref{pathc}}. 
Let $\fq$ be the Riemann map of $\om$ normalized as in \rife{fq2}.
By Remark \ref{scale} we can assume without loss of generality that $a_1=\fq^{'}(0)=1$.
Obviously $f(tz)$ is univalent in $\mb{D}$ for any $t\in (0,1]$ and so
is $\frac{f(tz)-q}{t}$, which takes the form,
\beq\label{series1}
\frac{f(tz)-q}{t}=z+\sum\limits_{n=3}^{+\ii}a_n t^{n-1}z^n.
\eeq
Since $\pa \om_t=f(\{|z|=t\})$, $t<1$, then $\om_t$ is analytic for any $t\in (0,1)$. Also
$$
D_{\om_t}(0)=\sum\limits_{n=3}^{+\ii}\frac{n^2}{n-2}|a_n|^2t^2-|a_1|^2\leq  \sum\limits_{n=3}^{+\ii}\frac{n^2}{n-2}|a_n|^2-|a_1|=
D_{\om}(0)\leq 1,
$$
which by Theorem A shows that $\om_t\in \mathcal{A}_I$ for any $t$. Since the domain $\om$ is at least regular,
then by Remark \ref{rempom} $\fq$ admits a continuous extension to $\ov{\mb{D}}$.
Thus the series in \rife{series1} converges for $t=\pm 1$ and
$z\in \pa \mb{D}$, and then, for fixed $t\in (-1,1)$, it is totally convergent in $\{|z|\leq r\}$ for any $r\in (0,1)$ and
for fixed $z\in \mb{D}$ it is totally convergent in $\{|t|<\dt\}$, for any $\dt\in (0,1)$. Therefore $f(z,t)$ is a separately analytic function in
$\mb{D}\times (-1,1)$. By the Abel theorem the series converges uniformly in
$\ov{\mb{D}}\times [-1,1]$ to a continuous function. Then $f(z,t)$ is continuous in $\ov{\mb{D}}\times [-1,1]$.
However it is well known that a separately analytic and jointly continuous function is jointly analytic {(see e.g. \cite{Hor} Theorem 2.2.1)}, whence $f(z,t)$
is a jointly analytic in
$\mb{D}\times (-1,1)$. The first part of the claim readily follows since $f(z,0)=z$ and $f(z,1)=\fq(z)$.\\\\

Next, let $\mathcal{A}_{1}$ denote the topological space of $C^{1}$ simply connected domains endowed with the $C^1$-topology, with $\mathcal{A}_{1,I}$ the subset of domains of first kind and with $\mathcal{H}(\mb{D})$ the space
of Holomorphic funcions in $\mb{D}$. Then, let us define the map $F:\mathcal{H}(\mb{D})\times [0,1] \to \mathcal{H}(\mb{D})$ as follows,
$$
F(f,t)=f(\cdot,t)
$$
with $f(z,t)$ as in \rife{fqt} above. The induced map $\mathcal{F}:\mathcal{A}_1 \to \mathcal{A}_1$ takes the form,

$$
\mathcal{F}(\om,t)=\om_t=f(\mb{D},t),
$$
and obviously we have
$$
\mathcal{F}(\om,0)=\om_0=\mb{D} \quad \mbox{and} \quad \mathcal{F}(\om,1)=\om_1=\om.
$$

Therefore $\mathcal{I}:=\mathcal{F}(\cdot,1)$ is the identity map in $\mathcal{A}_1$ and in particular, by the first part of the claim,
$\mathcal{F}(\om,t)\in \mathcal{A}_{1,I}$ for any $t\in [0,1]$. We claim that the restriction
$\mathcal{F}:\mathcal{A}_{1,I}\times [0,1]\to\mathcal{A}_{1,I}$ is a deformation retract of the identity in the given topology,
which proves that $\mathcal{A}_{1,I}$ is contractible. By Remark \ref{rempom} we see that $f^{'}$ admits a continuous extension on
$\ov{\mb{D}}$ with $f^{'}(z)\neq 0$ on $\ov{\mb{D}}$. Since obviously $d_1(\om_0,\mb{D})=0$ and $d_1(\om_1,\om)=0$,
then to establish the claim it will be enough to prove that $d_1(\om_t,\om)$ is continuous in $[0,1]$.  We will prove a
statement which easily implies the claim, that is 
$$
d_{\ii}(t):=\sup\limits_{z\in \ov{\mb{D}}}|f(z,t)-z|+\sup\limits_{z\in \ov{\mb{D}}}|f^{'}(z,t)-1|
$$
is continuous in $[0,1]$. We argue by contradiction and assume that there exists $t_0\in [0,1]$ and sequences $t_{n,i}\to t_0\in [0,1]$,
$i=1,2$ such that $|d_{\ii}(t_{n,2})-d_{\ii}(t_{n,1})|\geq \eps_0$, for some $\eps_0>0$ and any $n$. { Clearly we can find sequences
$z_{n,i}$ and $w_{n,i}$, $i=1,2$ which are maximizers of the corresponding absolute values, such that for any $n$,
$$
d_{\ii}(t_{n,i})=|f(z_{n,i},t_{n,i})-z_{n,i}|+|f^{'}(w_{n,i},t_{n,i})-1|,\;i=1,2.
$$
Passing to suitable subsequences we can assume w.l.o.g. that $z_{n,i}\to z_{i,0}$ and $w_{n,i}\to w_{i,0}$, i=1,2 where obviously
$z_{i,0}$ and $w_{i,0}$ are maximizers of the corresponding absolute values for $t=t_0$. Consequently, as $n\to +\ii$,
we would find that
$$
\eps_0\leq |d_{\ii}(t_{n,2})-d_{\ii}(t_{n,1})|\leq
$$
$$
\left| |f(z_{n,2},t_{n,2})-z_{n,2}|+|f^{'}(w_{n,2},t_{n,2})-1|-|f(z_{n,1},t_{n,1})-z_{n,1}|-|f^{'}(w_{n,1},t_{n,1})-1|\right|\to
$$
$$
\left| |f(z_{2},t_0)-z_{2}|+|f^{'}(w_{2},t_0)-1|-|f(z_{2},t_0)-z_{2}|-|f^{'}(w_{2},t_0)-1|\right|=
$$
$$
|d_{\ii}(t_{0})-d_{\ii}(t_{0})|=0,
$$
which yields the desired contradiction.}

\bigskip

Finally let $\mathcal{A}_{0,1}$ denote the set of simply connected and regular domains with metric $d_{0,1}$,
with $\mathcal{A}_{0,1,I}$ the subset of domains of first kind and with
$\Gamma:\mathbb{S}^1\to \mathcal{A}_{0,1,I}$
be any continuous loop of the form $\Gamma(s)=f(z;s)$, $z\in \ov{\mathbb{D}}$, $s\in\mathbb{S}^1$,  where
each $f(\cdot;s)$ is normalized as in \rife{fq2} with $q_s=f(0;s)$. Then, by the first part of the statement,
the map $F:\mathcal{H}(\mb{D})\times [0,1]\times \mathbb{S}^1\mapsto \mathcal{H}(\mb{D})$,
$$
F(\cdot,t,s)=tq_s+\frac{f(\cdot,t;s )-q_s}{t},\quad (t,s)\in [0,1]\times \mathbb{S}^1,
$$
induces in $\mathcal{A}_{0,1,I}$ a continuous deformation of $\Gamma(\mathbb{S}^1)=\mathcal{F}({\mathbb{D}},1;\mathbb{S}^1)$ to
${\mathbb{D}}=\mathcal{F}({\mathbb{D}},0;\mathbb{S}^1)$. In other words any loop $\Gamma(\mathbb{S}^1)$ in $\mathcal{A}_{0,1,I}$
can be deformed continuously to ${\mathbb{D}}$, which shows
that $\mathcal{A}_{0,1,I}$ is simply connected. \finedim

\bigskip

\section{Appendix I} In this appendix we prove the Claim used in the proof of Theorem \ref{thm2}$(ii)$, that is,
if a sequence $\om_k$ of regular and simply connected domains satisfies
$d_{0,1}(\om,\om_k)\to 0$, where $\om$ is regular and simply connected, then  $\ov{\om}_k$ is uniformly locally connected.\\
We argue by contradiction and suppose that $\ov{\om}_k$ is not uniformly locally connected.
Then $\exists\,\eps_0>0$ such that
$\exists a_k,b_k\in \ov{\om}_k$ such that $|a_k-b_k|<\frac{1}{k}$ and $\{a_k,b_k\}\nsubseteq B_k$ for any compact and connected subset
$B_k\subset \ov{\om}_k$ such that diam$(B_k)<\eps_0$. Since $\om_k$ is uniformly bounded, then passing to a subsequence if necessary,
we can assume without loss of generality that there exists $z_0\in \mathbb{C}$ such that $a_k\to z_0$ and $b_k\to z_0$. By the kernel
convergence and since $a_k,b_k\in \ov{\om}_k$, then $z_0\in \ov{\om}$ and we are left with two possibilities: either
$z_0\in\om$ or $z_0\in\pa\om$. We can easily exclude the first case, since then any closed disk $\ov{B_r}(z_0)$ with $r$ small enough will
contain both $a_k$ and $b_k$ and satisfy $\ov{B_r}(z_0)\subset \om_k$, for any $k$ large enough, which is a contradiction. \\
If $z_0\in\pa\om$, since the domain is regular, then
we have two possibilities: either the boundary is locally $C^2$ near $z_0$ or $z_0=Q$, where $Q$ is one of the vertex points on $\pa\om$.
We discuss only the second case. The proof of the regular case follows exactly by the same argument but it is easier and
we omit it here to avoid repetitions. After suitable translations we can assume that $z_0=Q=0\in \Gamma_{1}\cap \Gamma_{2}$ with $\theta$
the inner angle of $\Gamma_1$ and $\Gamma_2$  at $0$ and where $\Gamma_j$ are the $C^2$ connected components of $\pa\om$ near $0$.
Since $\om$ is regular then we can find an univalent map $f:B_{\dt}(0)\cap \ov{\om}\to \mathbb{C}$  such that $B_{\dt}(0)\cap\pa\om$ is mapped to
a $C^2$ curve.\\
Taking a smaller $\dt$ if necessary and composing with a suitable univalent map, we can assume without loss of generality
that $f(0)=0$, $f(B_{\dt}(0)\cap\pa\om)=\{w\in \mathbb{C}\,|\,w\in (-1,1)\}$ and
$f(B_{\dt}(0)\cap {\om})\subset \{w\in \mathbb{C}\,|\,\mbox{arg}(w)\in (0,\pi)\}$, where $\om_\dt=B_{\dt}(0)\cap {\om}$ is a simply
connected set. Since $a_k\to 0$, then
$a_k\in B_\dt(0)$ for $k$ large enough, and since by assumption $f$ is continuous and univalent in $B_{\dt}(0)$, then $f(a_k)\to f(0)=0$.
Clearly the same holds for $b_k\to 0$ and so $f(b_k)\to 0$ and
for any $r$ there exists $\nu_r$ such that if $k>\nu_r$ then $f(a_k),f(b_k)\in \ov{B_{r}}(0)$.
In particular, since each $\om_k$ is regular and converges to $\om$ in the kernel sense,
then we can choose $r_0$ small enough such that $U_{k,r}=f^{-1}(\ov{B_{r}}(0))\cap \ov{\om}_k$ is connected for any $r\leq r_0$.
Since each $U_{k,r}$ is a compact subset of $\ov{\om}_k$ and diam$(U_{k,r})\to 0$  as
$r\to 0$, we conclude that $a_k,b_k$ are contained in a connected compact subset of $\ov{\om}_k$,
whose diameter is smaller than $\eps_0$ for any $k$, which is the
desired contradiction. \finedim

\bigskip

\section{Appendix II}
In this Appendix we prove \rife{A.a1}, \rife{A.b1}, \rife{A.c1}, that is, for any fixed $y\in \om$ we have
\beq\label{A.a}
R_{k}(x,y)\to R_{\ssb}(x,y) \mbox{ in } C^{3}_{\rm loc}({\om})
\eeq
and,
\beq\label{A.b}
\gamma_k\to \gamma_{\ssb} \mbox{ in } C^{2}_{\rm loc}(\om),
\eeq
and moreover $h_k$ satisfies the uniform bound,
\beq\label{A.c}
|h_k(x)|\leq C, \; x\in B_{2\dt}(\pa \om_k)\cap \om_k.
\eeq
$\,$\\\\
{\em Proof of } \rife{A.a}.
We can assume without loss of generality that $y=0$. Clearly,  since by assumption $\om_k\to \om$ in the $C^1$-topology, we have that
$0\in\om_k$ for any $k$ large enough. Let $G_k(x,y)$ be the Green's function for $\om_k$ and $R_k(x,y)$ its regular part.
Since $d_1(\om_k,\om)\to 0$ then for any $\dt>0$ we have that $\pa\om_k\subset B_\dt(\pa\om)$ for any $k$ large and in particular
there exists $t_k\to 0^+$, and, for each $k$, $\eps_k>0$ and a one-to-one map $\Phi_k:B_{2\dt}(\ov{\om})\times [0,1]\to B_{\dt}(\ov{\om})$, that  satisfy
$$
\Phi_k(x,t)=x+tV_k(x)+\mbox{o}(t),\, t\in[0,t_k+\eps_k),
$$
$$
\Phi_k\in C^{1}(B_{2\dt}(\ov{\om})\times [0,1];\;B_{\dt}(\ov{\om})),\quad \|V_k\|_{C^{1}(B_{2\dt}(\ov{\om});\;B_{\dt}(\ov{\om}))}\leq C,
$$

$$
\lim\limits_{t\to 0^+}\sup\limits_{x\in B_{2\dt}(\ov{\om})}\frac{|\mbox{o}(t)|}{t}=0,
$$

and
$$
\Phi_k(\om,0)=\om,\qquad\Phi_k(\om,t_k)=\om_k.
$$
Let $\om_k(t)=\Phi_k(\om,t)$, $t\in [0,t_k+\eps_k)$ and $G_{\sscp \om_k(t)}(x,y)$ be the corresponding Green's function.
At this point we can apply a result in \cite{H} (Example 3.4) which shows that the map $t\mapsto G_{\sscp \om_k(t)}(x,y)$ is
differentiable for $x\neq y$ and in particular that the Hadamard variational formula holds,
$$
\frac{\pa }{\pa t}G_{\sscp \om_k(t)}(x,y)=-\int\limits_{\pa \om_k(t)}\frac{\pa G_{\sscp \om_k(t)}}{\pa \nu_{z}}(x,z)
\frac{\pa G_{\sscp \om_k(t)}}{\pa \nu_{z}}(z,y)
<V_k(z),\nu_z>d\sg(z),\;t\in[0,t_k+\eps_k),
$$
for $\{ x,y \}\in \om_k(t)$, see also \cite{Pe}. Here $\nu_z$ is the unit outer normal to $\om_k(t)$. Let $B_R=B_R(z)$ be any relatively compact disk
in $\om$ such that $0\notin \ov{B_R}$.
Clearly for $k$ large we have
$\ov{B_R}\cup \{0\}\subset{\om_k(t)}$, $\forall\,t \in[0,t_k+\eps_k)$ and then we can write,
$$
G_{k}(x,0)=G_{\ssb}(x,0)+\left(\frac{\pa }{\pa t}G_{\sscp \om_k(t)}(x,0)\right)_{t=0}t_k+\mbox{o}_x(t_k),\quad x\in \ov{B_R},
$$
where
$$
\left(\frac{\pa }{\pa t}G_{\sscp \om_k(t)}(x,0)\right)_{t=0}=-\int\limits_{\pa \om}\frac{\pa G_{\ssb}}{\pa \nu_{z}}(x,z)
\frac{\pa G_{\ssb}}{\pa \nu_{z}}(z,0)
<V_k(z),\nu_z>d\sg(z),
$$
and $\mbox{o}_x(t_k)$ is an infinitesimal quantity which satisfies
\beq\label{doublecheck}
\forall x\in \ov{B_R}, \quad \lim\limits_{k\to+\ii}\frac{|\mbox{o}_x(t_k)|}{t_k}=0.
\eeq

Since $\om$ is of class  $C^1$, then,
$$
\sup\limits_{z\in\pa\om}\left|\frac{\pa G_{\ssb}}{\pa \nu_{z}}(x,z) \frac{\pa G_{\ssb}}{\pa \nu_{z}}(z,0)<V_k(z),\nu_z>\right|\leq
{C_R},
$$
where ${C_R}$ depends only by $R$ and $\om$. As a consequence we conclude in particular that for any $ x\in\ov{B_R}$
it holds,

$$
\left|G_{k}(x,0)-G_{\ssb}(x,0)\right|\leq C_R t_k(1+o_x(1))\to 0, \mbox{ as }k\to 0.
$$

At this point we observe that, for  any smooth domain $\om_1$ lying in the interior of $\om$ and satisfying $0\notin \pa\om_1$,
we have that $R_k(x,0)$ is the unique solution of
$$
\left\{
\begin{array}{lll}
-\Delta R_{k}(x,0)&=& 0\quad \mbox{in}\quad \om_1, \\
\hspace{0.3cm}R_{k}(x,0)&=& G_{k}(x,0)+\frac{1}{2\pi}\log(|x|) \quad \mbox{on}\quad \partial \om_1,
\end{array}\right.
$$

and $R_{\ssb}(x,0)$ is the unique solution of
$$
\left\{
\begin{array}{lll}
-\Delta R_{\ssb}(x,0)&=& 0\quad \mbox{in}\quad \om_1, \\
\hspace{0.3cm}R_{\ssb}(x,0)&=& G_{\ssb}(x,0)+\frac{1}{2\pi}\log(|x|) \quad \mbox{on}\quad \partial \om_1.
\end{array}\right.
$$
Since $\pa \om_1$ is compact,  it can be covered with a finite number of balls $B_{R_j}(z_j)$, $z_j\in \pa\om$, $j=1,\cdots,N$,
such that
$0\notin \ov{B_{R_j}(z_j)}$. As a consequence $R_{k}(x,0)-R_{\ssb}(x,0)$ is harmonic in $\om_1$ and
$|R_{k}(x,0)-R_{\ssb}(x,0)|\leq C_{1,x} t_k\to 0$, as $k\to +\ii$, where $C_{1,x}=\max\limits_{j}\{C_{R_j}\}(1+\mbox{o}_x(1))$. Therefore $R_{k}(x,0)-R_{\ssb}(x,0)$
converges to $0$ pointwise on $\pa\om_1$ and then also in $C^m_{\rm loc}(\om_1)$ for any $m\geq 1$.
Since $\om_1$ is arbitrary, then the proof of \rife{A.a} is completed.\\\\

{\em Proof of} \rife{A.b}.
We first observe that actually,
\beq\label{A.a.1}
R_{k}(x,y)\to R_{\ssb}(x,y) \mbox{ in } C^{3}_{\rm loc}({\om}\times {\om}).
\eeq

Indeed from \rife{A.a}, and since $R_k(x,y)=R_k(y,x)$, then for fixed $x\in \om$, $R_k(x,y)\to R_{\ssb}(x,y)$ in $C^3_{\rm loc}(\om)$.
Then, since $R_k(x,y)$ is harmonic, it is not difficult to check that \rife{A.a.1} holds.\\
As a consequence, since $\gamma_k(x)-\gamma_{\ssb}(x)=R_k(x,x)-R_{\ssb}(x,x)=\lim\limits_{y\to x} (R_k(x,y)-R_{\ssb}(x,y))$,
then, passing to the limit as $k\to +\ii$, we see that because of \rife{A.a.1} we can actually exchange the limits,
to conclude that $\gamma_k(x)\to \gamma_{\ssb}(x)$ pointwise and in particular locally uniformly in $\om$. The same
argument works for the derivatives, since for example $\nabla(\gamma_k(x)-\gamma_{\ssb}(x))=2\nabla (R_k(x,x)-R_k(x,x))$, which
concludes the proof of \rife{A.b}.\\\\

{\em Proof of} \rife{A.c}.
Since $\inf\limits_{k}\,$dist$(0, \pa\om_k)\geq d>0$, by \rife{A.b} we are reduced to prove that, for any $\dt>0$ small enough, we have
$$
R_k(x,0)\leq C, \; x\in B_{2\dt}(\pa \om_k)\cap \om_k.
$$
However this is obvious since $R_k(x,0)$ is harmonic in $\om_k$ and satisfies $2\pi R_k(x,0)=-\log(|x|)$ for $x\in \pa \om_k$. Therefore
$2\pi R_k(x,0)\leq \sup\limits_{x\in \pa\om}(-\log(|x|))\leq -\log{d}$ for any
$x\in B_{2\dt}(\pa \om_k)\cap \om_k$, whenever $\dt>0$ satisfies $3\dt<d$. \finedim

\bigskip

\section{Appendix III}
We discuss an example which shows that \rife{nan} is not necessary for a domain $\om$ to be of first kind.
It is also shown that, increasing each $|a_n(t)|$ along certains path, one can get well inside $\mathcal{A}_{II}$.
Indeed, let us consider the following family of functions,
$$
f(z;t)=z+t \frac{2}{3}\frac{z^3}{3}+t \frac{1}{3}\frac{z^5}{5},\;|z|\leq 1, \; t\in \left[0,\frac{5}{2}\right],
$$
which satisfies,
$$
\sum\limits_{n=3}^{+\ii}n|a_n|=3|a_3(t)|+5|a_5(t)|=t,
$$
and
$$
D_{\ssb}(0)=\sum\limits_{n=3}^{+\ii}\frac{n^2}{n-2}|a_n(t)|^2-|a_1(t)|^2=\frac{13}{27}t^2-1.
$$
Some elementary numerics shows that for $t\in [0,\frac{5}{2}]$, $f(z,t)$ is univalent and maps $\ov{\mathbb{D}}$ onto a $C^1$
and symmetric (w.r.t. the $x$ and $y$ axis)
domain $\ov{\om_t}=f(\ov{\mathbb{D}},t)$ such that, putting $t_0=3\sqrt{\frac{3}{13}}\simeq 1,44$, it holds:\\
if $t\leq t_0$ then $\om_t\in\mathcal{A}_I$;\\
if $t\in\left(t_0,\frac{3}{2}\right]$ then  $\om_t\in\mathcal{A}_{II}$ but it is starlike and $\gamma_{\sscp \om_t}$ has a unique maximum point;\\
if $t\in\left(\frac{3}{2},\frac{5}{2}\right]$, $\om_t\in\mathcal{A}_{II}$, $\gamma_{\sscp \om_t}$ has two maximum points and, for
$t$ close enough to $\frac{5}{2}$, $\om_t$ is a dumbbell shaped non-starlike domain.\\\\

\end{document}